\documentclass[12pt]{amsart}
\usepackage{amscd,amssymb}
\usepackage[arrow,matrix]{xy}
\usepackage[plainpages,backref,urlcolor=blue]{hyperref}

\theoremstyle{plain}
\newtheorem{thm}[subsection]{Theorem}

\newtheorem{prop}[subsection]{Proposition}

\theoremstyle{definition}
\newtheorem{rk}[subsection]{Remark}

\newtheorem{ex}[subsection]{Example}

\numberwithin{equation}{section}
\newcommand{\A}{{\mathcal A}}

\newcommand{\R}{\mathbb{R}}
\newcommand{\C}{\mathbb{C}}

\newcommand{\PP}{\mathbb{P}}

\DeclareMathOperator{\codim}{codim}

\begin{document}

\title [ Monodromy of triple point line arrangements ]
{Monodromy of triple point line arrangements }

\author[Alexandru Dimca]{Alexandru Dimca$^1$}
\address{Institut Universitaire de France et Laboratoire J.A. Dieudonn\'e, UMR du CNRS 7351,
                 Universit\'e de Nice Sophia-Antipolis,
                 Parc Valrose,
                 06108 Nice Cedex 02,
                 France}
\email{dimca@unice.fr}

\thanks{$^1$ Partially supported by the  ANR-08-BLAN-0317-02 (SEDIGA)} 



\begin{abstract}
We show that the monodromy operator action on the first cohomology group of the Milnor fiber is combinatorially determined
for line arrangements with at most triple points and containing at most 18 lines, with one possible exception. This last case depends on the existence of a line arrangement involving 18 lines and 48 triple points with a special distribution with respect to smooth conics.

\end{abstract}

\maketitle


\section{Introduction} \label{sec1}

Let $\A$ be an arrangement of $d$ lines in $\PP^{2}$, with $d \geq 2$, given by a reduced equation 
$ Q(x)=0$. 
Consider the corresponding complement $M$ defined by $Q(x)\ne 0$ in $\PP^{2}$, and the global Milnor fiber $F$ defined by $Q(x)-1=0$ in $\C^3$ with monodromy action $h:F \to F$, $h(x)=\exp(2\pi i/d)\cdot x$.

\medskip

It is an open question whether the monodromy operator $h^1:H^1(F) \to H^1(F)$ is combinatorially determined, i.e. determined by the intersection lattice $L(\A)$. The $1$-eigenspace $H^1(F)_1$ coincide with $H^1(M)$, and hence it is known to be determined by the lattice $L(\A)$, see \cite{OT}. The answer is not known even for the first Betti number $b_1(F)$. Several interesting examples have been computed by D. Cohen and A. Suciu, \cite{CS}, \cite{S1}.

\medskip
We assume in the sequel that the line arrangement $\A$ has only double and triple points. Then it is known that $h^1:H^1(F) \to H^1(F)$ is trivial unless $d=3m$ for some integer $m\geq 1$, and then only the eigenvalues $1$, $\epsilon= \exp(2\pi i/3)$ and $\epsilon ^2$ are possible, see for instance \cite{D2}, Cor. 6.14.15. This result gives also an upper bound on the multiplicity of $\epsilon$ as an eigenvalue of $h^1$, namely
$(3m-1)/2$, which is very poor, as the examples below show. For other related upper bounds see \cite{CL}.

Moreover, in a recent paper \cite{L1}, A. Libgober has followed the approach started in \cite{CL} and has shown that if $h^1 \ne Id$, then necessarily 
$\A$ is composed of a reduced pencil, i.e. there is a pencil $aQ_1+bQ_2$ of curves of degree $m$ in $\PP^2$
such that if we set $Q_3=Q_1+Q_2$, then $Q=Q_1\cdot Q_2\cdot Q_3$. The converse implication, that a line arrangement composed of a reduced pencil has a non-trivial $h^1$ is easy, see for instance Theorem 3.1 in \cite{DP}.

Set $\A_i:Q_i=0$ for $i=1,2,3$ and note that $\A=\A_1\cup \A_2 \cup \A_3$.
Let $T_0$ be the base point set of the pencil, i.e. the $m^2$ points of intersection of the curves $C_1$ and $C_2$, given by the union of all lines in $\A_1$ (resp. in $\A_2$). Let $T_i$ be the set of triple points of the arrangement $\A_i$, for $i=1,2,3$. Then $T=T_0 \cup T_1 \cup T_2 \cup T_3$
is the set of triple points of the arrangement $\A$.

Our result is the following.

\begin{thm} \label{Thm1} 
Let $\A$ be an arrangement of $d$ lines in $\PP^{2}$, with $d=3m$, such that $\A$ has only double and triple points and $\A$ is composed of a reduced pencil. Then, the following hold.

\noindent (i) If $m<6$ or if $m=6$ and $|T|<48$, then the
monodromy operator $h^1:H^1(F) \to H^1(F)$ is combinatorially determined. Moreover, one has
\begin{equation} \label{eq1}
\dim H^{1,0}(F)_{\epsilon}=\dim H^{0,1}(F)_{\epsilon^2}=1
\end{equation}
except for the case of a Ceva type arrangement ($m=3$, $|T_1|=|T_2|=|T_3|=1$), when
\begin{equation} \label{eq2}
\dim H^{1,0}(F)_{\epsilon}=\dim H^{0,1}(F)_{\epsilon^2}=2.
\end{equation}

\noindent (ii) If $m=6$ and $|T|=48$, which is the maximal possible value for $|T|$,
then each of the  subarrangements $\A_i$ has a set of $4$ triple points $T_i$.
When each set of the three sets  of $8$ points given by $T_1\cup T_2$, $T_2\cup T_3$ and $T_1\cup T_3$ is situated on a conic (necessarily smooth), then
\begin{equation} \label{eq3}
\dim H^{1,0}(F)_{\epsilon}=\dim H^{0,1}(F)_{\epsilon^2}=2.
\end{equation}
Otherwise, one has
\begin{equation} \label{eq4}
\dim H^{1,0}(F)_{\epsilon}=\dim H^{0,1}(F)_{\epsilon^2}=1.
\end{equation}

\end{thm}

It is known that $H^1(F)_{\ne 1}$ is a pure Hodge structure of weight 1, see \cite{BDS} and \cite{DP}
for two distinct proofs, and hence the notation $H^{1,0}(F)_{\epsilon}$ is unambiguous.

\begin{rk}
\label{rkexamples} 
An exceptional  arrangements $\A$ with 18 lines as described in the second part of the above Theorem  and satisfying \eqref{eq3} does exist, see Example \ref{keyex} below.
On the other hand,
we do not know whether exceptional  arrangements $\A$ with 18 lines satisfying \eqref{eq4} really exist. In the affirmative case they would provide a counter example to the determinancy of $b_1(F)$ or $h^1$ by the lattice $L(\A)$.

\medskip

In fact, the only arrangements having at most triple points and composed of a reduced pencil that were commonly known, before Example \ref{keyex} was kindly provided by  Masahiko Yoshinaga, are the following.

\medskip

\noindent(i) ($m=1$) three concurrent lines, where we may take $Q_1=x-y$, $Q_2=y-z$.

\medskip

\noindent(ii) ($m=2$) the $A_3$-arrangement, where we may take $Q_1=x^2-y^2$, $Q_2=y^2-z^2$. Note that each of the subarrangements $\A_i$ in Theorem \ref{Thm1}, (ii) is linearly equivalent to this $A_3$-arrangement.

\medskip

\noindent(iii) ($m=3$) the Ceva arrangement, where we may take $Q_1=x^3-y^3$, $Q_2=y^3-z^3$.

\medskip

\noindent(iv) ($m=3$) the Hesse arrangement, where we take $\A$ to be the union of three singular fibers out of the four singular fibers of the Hesse pencil $a(x^3+y^3+z^3)+bxyz$.

\medskip

\noindent(v) ($m=4$) the generic hyperplane section of  the $D_4$-plane arrangement in $\PP^3$ given by the equation
$$(x^2-y^2)(x^2-z^2)(x^2-t^2)(y^2-z^2)(y^2-t^2)(z^2-t^2)=0.$$
The pencil structure in this case comes from the identity
$$(a-b)(c-d)+(a-c)(d-b)+(a-d)(b-c)=0,$$
where we set $a=x^2$, $b=y^2$, $c=z^2$ and $d=t^2$. The values given above in \eqref{eq1} were obtained in the case $D_4$ in a different way by A. M\u acinic and S. Papadima in \cite{MP}.

\medskip

\end{rk}

\begin{rk}
\label{rkcounterexamples} 
Note that for arrangements having points of multiplicity  $4$, Libgober's result is no longer true: for instance the arrangement 
$$xyz(x^4-y^4)(y^4-z^4)(x^4-z^4)=0$$
satisfies $h^1 \ne Id$ (coming from the existence of the non-reduced pencil $Q_1=x^4(y^4-z^4)$ and $Q_2=y^4(z^4-x^4)$ see Remark 3.4. (iii) in \cite{DP} for $r=4$), but it is not composed of a reduced pencil. 
\end{rk}
 
\medskip

We would like to thank  Masahiko Yoshinaga for providing the key Example \ref{keyex}.

\section{Proof of Theorem \ref{Thm1}} \label{sec2}

\proof

The case $m=1$ is obvious, so we assume $m \geq 2$ in the sequel.

Let $S=\C[x,y,z]$ be the graded ring of polynomials in $x,y,z$. Consider now the evaluation map
\begin{equation} \label{eval}
\rho: S_{2m-3} \to \C^T
\end{equation}
obtained by picking up a representative $s_t$ in $\C^3$ for each point $t\in T$ and sending
a homogeneous polynomial $h \in S_{2m-3}$ to the family $(h(s_t))_{t\in T}$.

For any line arrangement $\A$ with at most triple points, it follows from Theorem 2 in \cite{BDS} (and the discussion folowing it) that 
\begin{equation} \label{eval2}
\dim H^{1,0}(F)_{\epsilon}=\dim H^{0,1}(F)_{\epsilon^2}=\dim (Coker \rho),
\end{equation}
while
\begin{equation} \label{eval3}
\dim H^{1,0}(F)_{\epsilon^2}=\dim H^{0,1}(F)_{\epsilon}=0.
\end{equation}
Alternatively, one may use the approach described in \cite{D0}, Chap.6., Thm. 4.15.

Note that
$$\dim (Coker \rho)=|T|- \codim (ker \rho)$$
which by definition is the {\it superabundance} $s_{2m-3}(T)$ of the finite set of points $T$ with respect to the polynomials in $S_{2m-3}$.

Hence, to determine the monodromy operator $h^1:H^1(F) \to H^1(F)$ (and the corresponding mixed Hodge structure) boils down to computing the dimension of the kernel of $\rho$. 

Let $h \in S_{2m-3}$ satisfy $\rho(h)=0$. In particular, it follows that $h$ vanishes at the points in $T_0$ and a direct application of Noether's $AF+BG$ Theorem, implies that 
\begin{equation} \label{rel1}
h=h_1Q_1+h_2Q_2,
\end{equation}
for some polynomials $h_1,h_2 \in S_{m-3}$. Moreover, such a pair $(h_1,h_2)$ is unique: an equality
$h_1Q_1+h_2Q_2=k_1Q_1+k_2Q_2$ yields $Q_2$ divides $h_1-k_1$, and hence looking at the degrees we see that $h_1=k_1$.

If $m=2$, this means the kernel of $\rho$ is trivial, so the claim follows, in particular in this case (which corresponds to the $A_3$-arrangement in Remark \ref{rkexamples}), we get 
\begin{equation} \label{eval4}
\dim H^{1,0}(F)_{\epsilon}=\dim H^{0,1}(F)_{\epsilon^2}=\dim (Coker \rho)=1.
\end{equation} 

Assume from now on that $m\geq 3$. For $t_1 \in T_1$,  we get using \eqref{rel1} $$0=h(t_1)=h_2(t_1)Q_2(t_1),$$ 
and hence
$h_2(t_1)=0$. Indeed, $Q_2(t_1)\ne 0$, since there is no line in a subarrangement $\A_i$ passing through a triple point in a set $T_k$ for $k \ne i$. Similarly, for $t_2 \in T_2$ we get $h_1(t_2)=0.$

For $t_3\in T_3$ we get two relations: $0=Q_3(t_3)=Q_1(t_3)+Q_2(t_3)$ and $0=h(t_3)=h_1(t_3)Q_1(t_3)+h_2(t_3)Q_2(t_3)$. Since $Q_1(t_3)\ne 0$ and $Q_2(t_3)\ne 0$, it follows that $h_1(t_3)=h_2(t_3)$.

Hence $\dim Ker(\rho)=\dim H$, where
$$H=\{(h_1,h_2) \in S_{m-3}^2~~:~~h_1|T_2=0,~~h_2|T_1=0 \text{ and } (h_1-h_2)|T_3=0\}.$$

Now we have to see when the vanishing conditions in the definition of $H$ are linearly independent conditions. If this is so, then the dimension of $H$ is determined by the number of these conditions, which is $|T|-m^2$, so a combinatorially determined number. 
And then we have to show that the special cases when these conditions are not linearly independent
(and in fact the ranks of the corresponding linear systems) are also combinatorially determined.

This is done by a case-by-case analysis, based on the following approach. Consider the evaluation mapping
$$ \rho': S^2_{m-3} \to \C^{T_1 \cup T_2 \cup T_3}$$
given by sending a pair $(h_1,h_2) \in S^2_{m-3}$ to the families of complex numbers $$((h_1(t_2))_{t_2 \in T_2}, (h_2(t_1)_{t_1 \in T_1}, (h_1(t_3)-h_2(t_3))_{t_3 \in T_3}).$$

\begin{prop}
\label{prop1}
\noindent (i) $\dim (Coker \rho)= \dim (Coker \rho')+1.$

\medskip
 
\noindent (ii) $\dim (Coker \rho')=|T_1|+|T_2|+|T_3|-\codim H $ and  
$$\dim (Coker \rho') \leq \min\{s_{m-3}(T_1 \cup T_2), s_{m-3}(T_1 \cup T_3), s_{m-3}(T_2 \cup T_3)\}.$$
 
\end{prop}

\proof The first claim is a direct computation based on the equality
$$\dim S_{2m-3}+1=2\dim S_{m-3}.$$
The first equality in the second claim is just the definition of $\dim (Coker \rho')$, where the codimension of $H$ is computed with respect to the ambient space $S^2_{m-3}$.

Now we prove the inequality $\dim (Coker \rho') \leq s_{m-3}(T_1 \cup T_2)$, and note that the other two cases are completely symmetrical.

Set $V_1=\{h \in S_{m-3}~~:~~h|T_1=0\}$ and $V_2=\{h \in S_{m-3}~~:~~h|T_2=0\}$. Note that $s_{m-3}(T_j)=0$ for any $j=1,2,3$, see \cite{BDS}, the discussion following Theorem 2. These vanishings are related to the computation of the spectrum of a hyperplane arrangement in \cite{BS}.

In fact, $s_{m-3}(T_j)=0$ is equivalent to the surjectivity of the corresponding evalation mapping, and the vanishing $s_{a-3}(T_j)=0$ is claimed in \cite{BDS} for $2m/3<a <m$. However, it is clear that the surjectivity claim for $a$ implies the surjectivity claim for $a+1$, e.g. by choosing the coordinates such that all the points in $T_j$ are in the affine chart $z=1$ and using $\C[x,y]_{\le a} \subset \C[x,y]_{\le a+1}$. The special case $m=3$ follows by a direct inspection.

It follows that $\codim V_j= |T_j|$ for $j=1,2$ and hence 
$$\codim(V_1+V_2)=\codim V_1+ \codim V_2-\codim (V_1 \cap V_2)=s_{m-3}(T_1\cup T_2).$$
Consider now the restriction of the evaluation mapping $\rho'$ to the subspace $V_2 \times V_1 \subset S^2_{m-3}$, namely $\theta: V_2 \times V_1 \to \C^{T_3}$, $(h_1,h_2) \mapsto (h_1(t_3)-h_2(t_3))_{t_3 \in T_3}.$

Using $ker \rho'=ker \theta$, one has 
$$
\dim (Coker \rho')=|T_1|+|T_2|+|T_3|-\codim (ker \rho')= 
|T_3| - 2 \dim S_{m-3}+|T_1|+|T_2|+$$
$$+\dim (ker \theta)= |T_3|-\dim(V_2 \times V_1)+\dim (ker \theta)=\dim (Coker \theta).$$

On the other hand, the map $\theta$ factors through $V_1+V_2$. Since $s_{m-3}(T_3)=0$, it follows that
$\dim (Coker \theta) \leq \codim (V_1+V_2)=s_{m-3}(T_1\cup T_2)$.

\endproof


Now we discuss the various possible cases.

\bigskip

{\bf Case $m=3$}.

In this case, the subarrangements $\A_i$ are either given by a triangle (in which case $T_i=\emptyset$)
or by 3 concurrent lines (when $|T_i|=1$).

If there is at least one subarrangement of the first type, then we have at most two equations involving the two constants $h_1$ and $h_2$, hence the conditions are independent. In particular in this case (which covers  the Hesse arrangement in Remark \ref{rkexamples}, (iv)), we get 
\begin{equation} \label{eval5}
\dim H^{1,0}(F)_{\epsilon}=\dim H^{0,1}(F)_{\epsilon^2}=\dim (Coker \rho)=1.
\end{equation} 

On the other hand, if all subarrangements are of the second type (this corresponds exactly to the condition $|T|=9+3=12$), then we have three equations involving the two constants $h_1$ and $h_2$, hence the conditions are not independent, one equation being superfluous. In particular in this case (which covers  the Ceva arrangement in Remark \ref{rkexamples}, (iii)) we get 
\begin{equation} \label{eval6}
\dim H^{1,0}(F)_{\epsilon}=\dim H^{0,1}(F)_{\epsilon^2}=\dim (Coker \rho)=2.
\end{equation} 

\bigskip

{\bf Case $m=4$}.

In this case, the subarrangements $\A_i$ can be either generic, i.e. $T_i=\emptyset$, or special,
when $|T_i|=1$.  As two points impose independent conditions on linear forms,
 it follows that $s_1(T_i \cup T_j)=0$ for all $i,j$. Hence in any case we get
\begin{equation} \label{eval8}
\dim H^{1,0}(F)_{\epsilon}=\dim H^{0,1}(F)_{\epsilon^2}=\dim (Coker \rho)=1.
\end{equation}

\bigskip

{\bf Case $m=5$}.

In this case, the subarrangements $\A_i$ can be either generic, i.e. $T_i=\emptyset$,  special,
when $|T_i|=1$ or even 2-special, when $|T_i|=2$. Note that in the last case, the line determined by the two triple points is in the arrangement $\A_i$. It follows that if two of the subarrangements are 2-special, then the corresponding family of 4 triple points are not collinear.

Now, at most 3 points or 4 non-collinear points impose independent conditions on conics. It follows that  $s_2(T_i \cup T_j)=0$ for all $i,j$, and hence again
\begin{equation} \label{eval9}
\dim H^{1,0}(F)_{\epsilon}=\dim H^{0,1}(F)_{\epsilon^2}=\dim (Coker \rho)=1.
\end{equation}

\bigskip

{\bf Case $m=6$}.

In this case, the subarrangements $\A_i$ can be either generic ($T_i=\emptyset$),  special
($|T_i|=1$), 2-special ($|T_i|=2$),  3-special ($|T_i|=3$)  or even 4-special ($|T_i|=4$).
Assume that the subarrangements $\A_i$ have been numbered such that the sum $\sigma=|T_1|+|T_2|$ is minimal among all the sums $|T_i|+|T_j|$ for $i\ne j$. 

If $\sigma \leq 4$, then the above argument  works again, since 4 points impose
independent conditions on cubics.

When $\sigma =5$, there are two possibilities (up-to a symmetry in $\A_1$, $\A_2$), namely
$(|T_1|,|T_2|)=(1,4)$ and $(|T_1|,|T_2|)=(2,3)$. In both cases, note that among these 5 points, there are at most 3 collinear ones (use the fact that for any 3-special or 4-special subarrangement of 6 lines, any line determined by 2 triple points belong to the subarrangement).
Then we may add two generic points to our 5 points in $T_1 \cup T_2$ and use the fact that 7 points impose independent conditions on cubics unless 5 are collinear, see \cite{GH}, p. 714.

When $\sigma =6$, there are two possibilities (up-to a symmetry in $\A_1$, $\A_2$), namely
$(|T_1|,|T_2|)=(2,4)$ and $(|T_1|,|T_2|)=(3,3)$. In both cases, note that among these 6 points, there are at most 3 collinear ones, so the same result applies.

When $\sigma =7$, there is just one possibility (up-to a symmetry in $\A_1$, $\A_2$), namely
$(|T_1|,|T_2|)=(3,4)$. Note that among these 7 points, there are at most 2 collinear ones, so the same result applies.

When $\sigma =8$, there is just one possibility, namely
$(|T_1|,|T_2|)=(4,4)$. Note that among these 8 points, there are at most 2 collinear ones. However, 8 points impose independent conditions on cubics unless 5 of them are collinear or all 8 points are situated on a conic, see \cite{GH}, p. 715. 

So, if at least one of the sets of 8 points given by $T_1\cup T_2$, $T_2\cup T_3$ and $T_1\cup T_3$ is not situated on a conic, then we conclude as before, since $V_1+V_2=S_3$.

In the opposite case, $H_0=V_1+V_2$ is a hyperplane in $S_3$, hence $\dim H_0=9$, and the evaluation map $\theta':H_0 \to \C^{T_3}$
is not surjective. Indeed, its kernel contains the cubic forms of the type $L_1q_{13}+L_2q_{23}$
where $L_1,L_2 \in S_1$, and $q_{13}$ (resp. $q_{23}$) is a fixed quadratic form vanishing on $T_1\cup T_3$ (resp. $T_2 \cup T_3$). This gives a 6-dimensional subspace in $H_0$, since the conics
$q_{13}$ and $q_{23}$ are not degenerated (by the above discussion on collinear points).
Hence in this case $\dim(Coker \theta)=\dim(Coker \theta')=1$, which yields the result.

\endproof

\begin{ex}
\label{keyex} (Yoshinaga's 18 line arrangement)
Let $c \in \R$ be a large real number, and set $a=\exp(2\pi i/6)$, a primitive $6$-th root of unity. Define the polynomial
$$Q_1(x,y,z)=(x^3-y^3)(x+y-cz)(ax+a^5y+cz)(a^5x+ay+cz).$$
A direct computation shows that
\begin{equation} \label{ex1}
Q_1(x,y,z)=x^6-y^6+3cx^4yz-3cxy^4z-c^3x^3z^3+c^3y^3z^3.
\end{equation} 
Consider now the cyclic permutation $\tau (x,y,z)=(y,z,x)$, acting in the obvious way on polynomials in $S$ and on the projective plane $\PP^2$.
Define new polynomials $Q_2(x,y,z)=\tau Q_1=Q_1(y,z,x)$ and $Q_3(x,y,z)=\tau^2Q_1=Q_1(z,x,y)$.

\medskip

Then we have the following.

\medskip

\noindent (i) Each of the three arrangements $\A_i:Q_i=0$ is equivalent to the $A_3$-arrangement. The corresponding sets of triple points are described as follows.
\begin{equation} \label{ex2}
T_1=\{(0:0:1), (-c:-c:1), (a^5c:ac:1), (ac:a^5c:1)\},
\end{equation} 
and $T_2=\tau T_1$, $T_3=\tau T_2$.

\medskip

\noindent (ii) Using the formula \eqref{ex1} we see that $Q_1+Q_2+Q_3=0$, i.e. we have a pencil
of degree 6. It is easy to check, using the asumption $c>>0$, that the base locus $T_0$ of this pencil is formed indeed by 36 points.

\medskip

\noindent (iii) The 8 points in $T_1 \cup T_2$ are situated on the smooth conic 
$$C_3:q_3(x,y,z)=y^2+cxz=0.$$
Similarly, the 8 points in $T_2 \cup T_3 =\tau(T_1 \cup T_2)$ are situated on the smooth conic 
$$C_1:q_1(x,y,z)=\tau q_3=z^2+cxy=0,$$
and the 8 points in $T_1 \cup T_3=\tau(T_2 \cup T_3)$ are situated on the smooth conic 
$$C_2:q_2(x,y,z)=\tau q_1=x^2+cyz=0.$$
It follows that for the line arrangement $\A:Q_1Q_2Q_3=0$ all the conditions implying the equalities in \eqref{eq3} hold.

\end{ex}

\end{document}